\documentclass[reqno,12pt,final]{amsart}

\usepackage{a4wide}
\usepackage{ucs}
\usepackage[utf8x]{inputenc}
\usepackage{amsmath,amssymb}

\usepackage[notref,notcite,color]{showkeys}
\usepackage{enumitem,color}

\newtheorem{theorem}                   {Theorem}

\newtheorem{lemma}           [theorem] {Lemma}

\newtheorem{fact}			 [theorem] {Fact}

\theoremstyle{remark}

\let\epsilon\varepsilon
\let\eps\varepsilon
\let\rho\varrho

\def\neighbor{\mathrm{N}}

\newcommand{\By}[2]{\overset{\mbox{\tiny{#1}}}{#2}}
\newcommand{\ByRef}[2]{   \By{\eqref{#1}}{#2} }

\newcommand{\gByRef}[1]{  \ByRef{#1}{>} }

\def\ex{\mathrm{ex}}

\newcommand{\order}[1]{v(#1)}

\newcommand{\rcolor}[1]{#1}
\newcommand{\bcolor}[1]{#1}

\title{The extremal function for partial bipartite
tilings}

\author{Codru\unichar{355} Grosu}
\address{Faculty of Automatic Control and Computers, Politehnica University of Bucharest, Splaiul Independen\unichar{355}ei 313, sector 6, Cod 060042, Bucharest, Romania}
\email{grosu.codrut@gmail.com}
\author{Jan Hladk\'y}
\address{Department of Applied Mathematics, Faculty of Mathematics and Physics, Charles University, Malostransk\'e n\'am\v est\'i
 25, 118 00, Prague, Czech Republic \and
 DIMAP, University of
 Warwick, Coventry, CV4 7AL, United Kingdom.}
  \email{honzahladky@gmail.com}
  \thanks{Most of the work was done while JH was
  affiliated with TU Munich. JH would like to thank the
  group of Anusch Taraz for nice
  environment. JH was partially supported by Grant Agency
  of Charles University, grant GAUK~202-10/258009, and by
  BAYHOST}

\begin{document}

\begin{abstract}
\PrerenderUnicode{\unichar{355}}
For a fixed bipartite graph $H$ and given $\alpha\in
(0,1)$ we determine the threshold $T_H(\alpha)$ which
guarantees that any $n$-vertex graph with at least
$T_H(\alpha){n \choose 2}$ edges contains
$(1-o(1))\frac\alpha{\order{H}}n$ vertex-disjoint copies of
$H$. \rcolor{In the proof we use a variant of a technique developed by Koml\'os~\bcolor{[Combinatorica 20 (2000),
203--218}]}.
\end{abstract}
\maketitle


\section{Introduction}
The Tur\'an Theorem~\cite{Tur}, 
\rcolor{one of the most important results} in Extremal Graph Theory, 
gives a sharp
threshold, denoted $\ex(n,K_r)$, for the maximum number of
edges of an $n$-vertex graph with no copy of $K_r$. Even
though the Tur\'an Theorem applies to any pair of values
$n$ and $r$, the interesting instances are rather those
when $n$ is large compared to $r$. Erd\H{o}s and
Stone~\cite{ErdosStone1946} extended the result by
determining the asymptotic behaviour of the
function $\ex(n,H)$ for a fixed non-bipartite graph $H$.
The same problem in the case that $H$ is a fixed
bipartite graph is
--- despite considerable effort --- wide open for most
graphs $H$. \bcolor{This is known as the \emph{Zarankiewicz problem}.} Let us recall that when $H$ has colour classes
of sizes $s$ and $t$, $s\le t$, then the
K\"ovari-S\'os-Tur\'an Theorem~\cite{KovariST} asserts
that
\begin{equation}\label{eq:KTS}
\ex(n,H)\le O(n^{2-1/s})=o(n^2)\;.
\end{equation}
 On the other hand, a standard random graph argument
 gives that $\ex(n,K_{s,t})\ge
 \Omega(n^{2-(s+t-2)/(st-1)})$.

It is natural to extend the above \emph{existential
questions} to \emph{tiling questions}. In such a~setting
one asks for the maximum number of edges of an $n$-vertex
graph which does not contain $\ell$ vertex-disjoint
copies of a graph $H$. This quantity \rcolor{is denoted by}
$\ex(n,\ell\times H)$. Erd\H os and
Gallai~\cite{ErdGall59} gave a complete solution to the
problem in the case when $H=K_2$.
\begin{theorem}[Erd\H os-Gallai, 1959]\label{thm:ErdGall}
Suppose that $\ell\le n/2$. Then
$$\ex(n,\ell\times
K_2)=\max\left\{(\ell-1)(n-\ell+1)+{\ell-1\choose
2},{2\ell-1\choose 2}\right\}\;.$$
\end{theorem}
Given $n,x\in\mathbb N$, $x\le n$, we define \rcolor{two graphs $M_{n,x}$ and 
$L_{n,x}$ as follows. 
The} graph
$M_{n,x}$ \rcolor{is} an $n$-vertex graph whose vertex set is
split into sets $A$ and $B$, $|A|=x, |B|=n-x$, $A$
induces a clique, $B$ induces an independent set, and
$M_{n,x}[A,B]\simeq K_{x,n-x}$. The graph $L_{n,x}$ is \rcolor{the complement of $M_{n,n-x}$,} \bcolor{i.e., it is} an
$n$-vertex graph whose edges induce a clique of order
$x$. Obviously,
$e(M_{n,\ell-1})=(\ell-1)(n-\ell+1)+{\ell-1\choose 2}$,
and $e(L_{n,2\ell-1})={2\ell-1\choose 2}$. Moreover, it is easy to check that there are \rcolor{no $\ell$ vertex-}disjoint
edges in either of the graphs $M_{n,\ell-1}$,
$L_{n,2\ell-1}$. Therefore, when $\ell<\frac25 n\rcolor{ + }O(1)$,
the graph $M_{n,\ell-1}$ is (the unique)
graph showing that $\ex(n,\ell\times
K_2)\ge (\ell-1)(n-\ell+1)+{\ell-1\choose
2}$. The graph
$L_{n,2\ell-1}$ is the unique extremal graph for the
problem otherwise.

Moon~\cite{Moon68Tiling} started the investigation
of $\ex(n,\ell\times K_r)$. Allen, B\"ottcher, Hladk\'y,
and Piguet~\cite{AllBottHlaPig} only recently determined
the behaviour of $\ex(n,\ell\times K_r)$ for the whole
range of $\ell$ in the case $r=3$, and they
made a substantial progress for larger values of
$r$. Simonovits~\cite{Simonovits99Turan} determined the value $\ex(n,\ell
\times H)$ for a non-bipartite graph $H$, fixed value of
$\ell$ and large~$n$. 

\rcolor{An equally important} density parameter which can be considered in the
context of tiling questions is the minimum degree of the host
graph. That is, we ask what is the largest possible
minimum degree of an $n$-vertex graph which does not
contain $\ell$ vertex-disjoint copies of $H$. In the case
$H=K_r$\rcolor{,} the precise answer  is given by
the Hajnal-Szemer\'edi Theorem\footnote{\bcolor{In its original formulation, the Hajnal-Szemer\'edi Theorem asserts that an $n$-vertex graph $G$ with minimum-degree at least $\frac{r-1}{r}n$ contains a $K_r$-tiling missing at most $r-1$ vertices of $G$, thus giving an answer only to the question of almost perfect tilings. When the minimum-degree of $G$ is lower, we can however add auxiliary vertices which are complete to $G$ and obtain an $n'$-vertex graph $G'$ such that the Hajnal-Szemer\'edi Theorem applies to $G'$. The restriction of the almost perfect $K_r$-tiling of $G'$ to $G$ gives a $K_r$-tiling which is optimal in the worst case.}}~\cite{hajnal70:_proof_p}.
An asymptotic threshold for a general fixed graph $H$ was determined by
Koml\'os~\cite{komlosTT}. \rcolor{In this case, the threshold depends} \bcolor{ on a parameter which Koml\'os calls the \emph{critical chromatic number}. The critical chromatic number of $H$ is a real between $\chi(H)-1$ and $\chi(H)$. Roughly speaking, graphs $H$ which possess a coloring with $\chi(H)$ colors with one of the color classes small, have the critical chromatic number close to $\chi(H)-1$. On the other hand, graphs $H$ which have only approximately balanced $\chi(H)$-colorings have the critical chromatic number close to $\chi(H)$. There is a natural way how to state our main result, Theorem~\ref{main_result}, using the critical chromatic number. However, we chose not to as in the bipartite setting of Theorem~\ref{main_result} it is possible to give a self-contained formula for the problem. Let us also note that Koml\'os' result~\cite{komlosTT} gives an asymptotic min-degree threshold even in the case when $H$ is bipartite. In this case the near-extremal graphs for the problem are complete bipartite graphs.}

In the present paper we \rcolor{use a variation of the technique developed by Koml\'os to} 
determine \rcolor{the} asymptotic
behaviour of the function $\ex(n,\ell\times H)$ for a
fixed bipartite graph~\bcolor{$H$}. Let $H$
be an arbitrary bipartite graph. Suppose that
$b:V(H)\rightarrow [2]$ is a proper coloring of $H$ which minimizes $|b^{-1}(1)|$. 
We define quantities
$s(H):=|b^{-1}(1)|$, $t(H):=|b^{-1}(2)|$. Obviously,
$s(H)\le t(H)$, and $s(H)+t(H)=\order{H}$. Furthermore, we define
$V_1(H):=b^{-1}(1)$ and $V_2(H):=b^{-1}(2)$. 
The sets
$V_1(H)$ and $V_2(H)$ are uniquely defined provided that
$H$ does not contain a balanced bipartite graph as \rcolor{one of its components; 
in this other case we fix a coloring $b$
satisfying the above conditions and use it to define uniquely $V_1(H)$ and $V_2(H)$.}

Given $s,t \in \mathbb N$, we
define a function $T_{s,t}:(0,1)\rightarrow (0,1)$ by
setting
\begin{equation}
\label{thresh_func}
	T_{s,t}(\alpha) := 
	\max\left\{ 
	\frac{2s\alpha}{s+t}\left(1 -
	\frac{s\alpha}{2(s+t)}\right),
	\alpha^2 \right\}\;,
\end{equation}
for $\alpha\in(0,1)$.
Note that $T_{s',t'}=T_{s,t}$ when $s'=ks$ and $t'=kt$.
Also, note that 
\begin{equation}\label{eq:TssVersusErdosGallai}
T_{s,s}(\alpha){n\choose
2}=\ex\left(n,\frac{\alpha n}2\times
K_2\right)+o(n^2)\;,
\end{equation} 
\rcolor{and, in general, \bcolor{for $s\le t$, the number}
$T_{s,t}(\alpha)\binom{n}{2}$ is asymptotically the maximum between 
the number of edges of $M_{n, \frac{\alpha s}{s+t}n}$ and $L_{n, \alpha n}$.}

Our main result is the following.
\begin{theorem}
\label{main_result}
Suppose that $H$ is a bipartite graph with no isolated
vertices, $s:=s(H), t:=t(H)$. Let $\alpha \in
(0,1)$ and $\eps > 0$. Then there exists an $n_0 = n_0(s,t,\alpha,\eps)$ 
such that for any $n \geq n_0$, any graph $G$ with 
$n$ vertices and at least $T_{s,t}(\alpha){n \choose 2}$ 
edges contains more than $(1- \eps)\tfrac{\alpha}{s+t}n$
vertex-disjoint copies of the graph $H$.
\end{theorem}

\rcolor{Let $H, s$ and $t$ be as in the hypothesis of the theorem,
$\epsilon' > 0$ and $\beta \in (0, 1)$. Then we may find an $\alpha > \beta(s+t)$ 
and an $\epsilon < \epsilon'$ sufficiently
small, such that for $n$ large enough, by Theorem~\ref{main_result}, 
any graph $G$ with $n$ vertices and at least 
$T_{s, t}(\alpha)\binom{n}{2} < (T_{s,t}(\beta(s+t))+\epsilon')\binom{n}{2}$ edges contains at least $\beta n$ vertex-disjoint
copies of $H$. Hence 
$\ex(n,\beta n\times H)\le T_{s,t}(\beta(s+t)){n\choose 2}+\eps' n^2$.}
This
asymptotically matches the lower bound which comes --- as
in Theorem~\ref{thm:ErdGall} --- from graphs $M_{n,\beta
\rcolor{s}n-1}$ and $L_{n,\beta\rcolor{(s+t)}n-1}$. 
\rcolor{Indeed, neither of these graphs contains $\beta n$ vertex-disjoint copies of $H$, as any such copy
would require at least $s$ vertices in the clique subgraph of $M_{n, \beta sn -1}$, and at least
$s+t=v(H)$ non-isolated vertices in $L_{n, \beta(s+t)n-1}$\bcolor{, respectively}.}
Note however
that for most values of \rcolor{$H$,} the graphs $M_{n,\beta
\rcolor{s}n-1}$ and $L_{n,\beta\rcolor{(s+t)}n-1}$ are not
extremal for the problem. For example, we can replace the
independent set in the graph $L_{n,\beta\rcolor{(s+t)}n-1}$ by any $H$-free graph. \bcolor{This links us to the Zarankiewicz problem}, and
suggests that an~\bcolor{exact} result is not within the reach of
current techniques.

The assumption on $H$ to contain no isolated vertices
in Theorem~\ref{main_result} is made just for the sake of
compactness of the statement. Indeed, let $H'$ be obtained
from $H$ by removing all the isolated vertices. Then there
is a 
simple relation \rcolor{between} the sizes of optimal coverings
by vertex disjoint copies of $H$ and $H'$ in an
$n$-vertex graph $G$. Let $x$ and $x'$ be the number
of vertices covered by a  maximum family of
vertex-disjoint copies of $H$ and $H'$ in $G$,
respectively. We have that
$$x=\min\left\{\order{H}\left\lfloor\frac{n}{\order{
H}}\right\rfloor,\frac{x'\order{H}}{\order{H'}}\right\}\;.$$

\smallskip

\bcolor{One can attempt to obtain an analogue of Theorem~\ref{main_result} for graphs with higher chromatic number. This however appears to be substantially more difficult. To indicate the difficulty, let us recall that there are two types ($M_{n,x}$ and $L_{n,x}$) of extremal graphs for the $H$-tiling \rcolor{problem} for bipartite $H$. The graphs $M_{n,x}$ and $L_{n,x}$ have a \emph{block structure}, i.e., 
\rcolor{their vertex set can be partitioned into \emph{blocks} (two, in this case), such that any two vertices from the same block have almost \bcolor{the same neighborhoods}. These two graphs}
 appear even in the simplest case \rcolor{of} $H=K_2$ (cf. Theorem~\ref{thm:ErdGall}). However, when $H$ is not balanced, \rcolor{if we let $\alpha$ go from $0$ to $1$, the transition between the two extremal structures which determine the threshold function occurs} at a different time in the evolution. 
On the other hand, there are five types of extremal graphs for the problem of determining $\ex(n,\ell\times K_3)$ as shown in~\cite{AllBottHlaPig}. All the five types have a block structure. It is plausible that when $H$ is a general $3$-colorable graph,
\rcolor{the same five types of extremal graphs determine the threshold function for $H$-tilings.}
However, the transitions between them occur at different times and the block \rcolor{sizes}
depend on various structural properties of $H$. In particular, we have indications that the critical chromatic number alone does not determine $\ex(n,\alpha n\times H)$ in this situation.
}

\bigskip

If $\mathcal{F}$ is a family of graphs, and $G$ is a graph,
an \emph{$\mathcal{F}$-tiling} in $G$ is a set of
vertex-disjoint subgraphs of $G$, each of them isomorphic to a graph in
$\mathcal{F}$. If $\mathcal{F}=\{H\}$ then we simply say
$H$-tiling. $V(F)$ denotes the vertices of $G$ covered by
an $\mathcal{F}$-tiling $F$, and $|F| = |V(F)|$ is the
\emph{size} of the tiling $F$. 
\rcolor{If $F$ is a collection of bipartite graphs,
we let $V_1(F) = \bigcup_{H \in F}V_1(H)$ and $V_2(F)=\bigcup_{H \in F}V_2(H)$.}
For $n \in \mathbb N$, we write $[n]$ to denote the set
$\{1,2,\ldots,n\}$.



\section{Tools for the proof of the main result}

Our main tool is Szemer\'edi's \rcolor{R}egularity \rcolor{L}emma (see
\cite{KS96,KuhnOsthusSurv} for surveys). To state it we
need some more notation.

Let $G = (V,E)$ be an $n$-vertex graph. If $A,B$ are
disjoint nonempty subsets of $V(G)$, the \textit{density}
of the pair $(A,B)$ is $d(A,B) = e(A,B)/(|A||B|)$. We say
that $(A,B)$ is an \emph{$\eps$-regular pair} if 
$
|d(X,Y) - d(A,B)| < \eps
$
for every $X
\subset A, |X| > \eps|A|$ and $Y \subset B, |Y| >
\eps|B|$.

The following statement asserts that large subgraphs of
regular pairs are also regular. 
\begin{lemma}
\label{lem:slicing}
Let $(A,B)$ be an $\eps$-regular pair with density $d$,
and let $A' \subset A, |A'| \geq
\alpha|A|, B' \subset B, |B'| \geq \alpha|B|$,
$\alpha\ge \eps$. Then $(A',B')$ is an $\eps'$-regular
pair with $\eps' = \max\{\eps/\alpha, 2\eps\}$, and for its density $d'$ we
have $|d'-d| < \eps$.
\end{lemma}

Let $\eps > 0$ and $d \in [0,1]$.  An $(\eps,
d)$-\textit{regular partition} of $G$  with
\textit{reduced graph} $R = (V',E')$ is  a partition $V_0
\dot\cup V_1 \dot\cup \ldots \dot\cup V_k$  of $V$ with
$|V_0| \leq \eps n$, $|V_i| = |V_j|$ for any  $1 \leq i <
j \leq k$, $V(R) = \{V_1, V_2, \ldots, V_k\}$,  such that
$(V_i,V_j)$ is an $\eps$-regular pair in $G$ of  density
greater than $d$ whenever $V_iV_j \in E(R)$, and  the
subgraph $G' \subset G$ induced by the $\eps$-regular 
pairs corresponding to the edges of $R$ has more than 
$e(G) - (d+3\eps)n^2/2$ edges. In this case, we also  say
that $G$ has an $(\eps, d)$-reduced graph $R$, and call 
the sets $V_i, 1 \leq i \leq k$, the \textit{clusters} of $G$.

The following lemma is a consequence of the so-called 
degree version of the Regularity Lemma \cite[Theorem 1.10]{KS96}.

\begin{lemma}[Regularity \rcolor{L}emma]
\label{lem:reg}
For every $\eps > 0$ and $m \in \mathbb N$ there is an $M=M(\eps, m)$ such that, if $G$ is any graph with more than $M$ vertices and $d \in [0,1]$ is any real number, then $G$ has an $(\eps, d)$-reduced graph $R$ on $k$ vertices, with $ m \leq k  \leq M$.
\end{lemma}

Given four positive numbers $a,b,x,y$ we say that the
pair $a,b$ \emph{dominates} the pair $x,y$, if
$\max\{x,y\}/\min\{x,y\}\ge \max\{a,b\}/\min\{a,b\}$. 
The following easy lemma states that $K_{a,b}$ has an
almost perfect $K_{s,t}$-tiling provided that $a,b$
dominates $s,t$.
\begin{lemma}
\label{lem:til}
For any $s,t\in\mathbb N$ there exists a constant $C$
such that the following holds. Suppose that the pair
$a,b\in\mathbb N$ dominates $s,t$. Then the graph
$K_{a,b}$ contains a $K_{s,t}$-tiling containing all but
at most $C$ vertices of $K_{a,b}$.
\end{lemma}
\begin{proof}
If $s=t$ then necessarily $a=b$. There
obviously exists a $K_{s,t}$-tiling containing all but at
most $C:=2(s-1)$ vertices of $K_{a,b}$.

With no loss of generality, we may suppose that
$a\le b$ and $s< t$. Then $as\le bt$ and $bs\le at$. A
tiling with $\lfloor (bt-as)/(t^2-s^2)\rfloor$ copies of
$K_{s,t}$ with the $s$-part of the $K_{s,t}$ placed in the
$a$-part of the $K_{a,b}$ and $\lfloor
(at-bs)/(t^2-s^2)\rfloor$ copies placed the other way
misses at most $C:=2(s+t-1)$ vertices of $K_{a,b}$.
\end{proof}

The next lemmas, versions of the Blow-up
Lemma~\cite{KSS_bl}, assert that regular pairs have
almost as good tiling properties as complete bipartite graphs.
\begin{lemma}
\label{lem:key}
For every $d>0,\gamma \in (0,1)$ and any two graphs $R$ and $H$, there is an $\eps = \eps(H,d,\gamma) >  0$ such that the
following holds for all positive integers $s$. Let $R_s$ be the graph obtained from $R$ by replacing every vertex of $R$ by $s$ vertices, and every edge of $R$ by a complete bipartite graph between the corresponding $s$-sets. Let $G$ be 
\rcolor{any} graph obtained 
similarly from $R$ by replacing \rcolor{every vertex of $R$ by $s$ vertices, and every edge of $R$ with an $\eps$-regular pair} 
of density at least $d$. If $R_s$ contains an $H$-tiling of
size at least $\gamma \order{R_s}$ then so does $G$.
\end{lemma}
\begin{lemma}
\label{lem:tiling2}
For every bipartite graph $H$ and every $\gamma,d>0$
there exists an $\eps = \eps(H,d,\gamma) >  0$ such that
the following holds. Suppose that there is an $H$-tiling in
$K_{a,b}$ of size $x$. Let $(A,B)$ be an arbitrary
$\eps$-regular pair with density at least $d$, $|A|=a$,
$|B|=b$. Then the pair $(A,B)$ contains an $H$-tiling of
size at least $x-\gamma (a+b)$.
\end{lemma}

Finally, let us state a straightforward corollary of the
K\"onig Matching Theorem.
\begin{fact}
\label{bip_matching}
Let $G = (A \dot\cup B, E)$ be a bipartite graph with color classes $A$ and $B$. If $G$ has no matching with $l+1$ edges, then $e(G) \leq l \max\{|A|, |B|\}$.
\end{fact}


\section{The proof}
In this section, we first state and prove the main
technical result, Lemma~\ref{main_lemma}. Then, we show how it
implies Theorem~\ref{main_result}.

For $s,t \in \mathbb N$, we set $\mathcal{F}_1 :=
\{K_{s,t},K_{s,t-1},K_2\}$ and $\mathcal{F}_2 :=
\{K_{st,t^2},K_{st-1,(t-1)t},K_{st,(t-1)t},K_2\}$. Let
us note that when $s < t$, the sizes of the two
color classes of any graph from $\mathcal{F}^* :=
\mathcal{F}_1 \cup \mathcal{F}_2$ dominate $s$ and $t$.

\rcolor{Let $F$ be a $K_{s,t}$-tiling in a graph $G$,
$s<t$. Suppose $E_0$ and $E_1$ are matchings in $G[V(G)-V(F), V_1(F)]$ and $G[V_2(F)]$, 
respectively, such that each copy $K$ of $K_{s, t}$ in $F$ has at most one vertex matched by $E_0$ and at most
one vertex matched by $E_1$. If any $K \in F$ which has a vertex matched by $E_0$, also has a vertex matched by $E_1$,
then we call the pair $(E_0, E_1)$ an \emph{$F$-augmentation}. Note that in this case $E_0$ and $E_1$ are vertex disjoint, as $V_1(F) \cap V_2(F)=\emptyset$.}

The main step in our proof of Theorem~\ref{main_result} is
the following lemma.
\begin{lemma}
\label{main_lemma}
Let $t > s \geq 1, \alpha \in (0,1)$ and $\eps > 0$. \rcolor{Then there exists an $\eps' = \eps'(s,t,\alpha,\eps) > 0$ and an $h = h(s,t,\alpha,\eps) > 0$ such that the following holds.} Suppose $G$ is an $n$-vertex graph with $n \geq h$ and $e(G) \geq T_{s,t}(\alpha){n \choose 2}$, and $F$ is a $K_{s,t}$-tiling in $G$ \rcolor{of maximum size with} $|F| \leq (1 - \eps)\alpha n$. Then one of the following is true:

\renewcommand{\theenumi}{\roman{enumi}}

\begin{enumerate}
\item there exists an $\mathcal{F}_1$-tiling $F'$ in $G$
with $|F'| \geq |F| + \eps'n$, or
\item 
there exists an $F$-augmentation $\rcolor{(}E_0, E_1\rcolor{)}$ such that
$E_0$ contains at least $\eps'n$ edges.
\end{enumerate}

\renewcommand{\theenumi}{\arabic{enumi}}
 
\end{lemma}
\begin{proof}
Set
\begin{equation*}
\label{eps_def}
\eps' :=\frac14 \min\left\{\frac{\eps \alpha^2}{3t+1},
\frac{\eps s \alpha}{(3t+1)(s+t)}\right\}\;,
\end{equation*}
and let $h$ be sufficiently large.

Suppose for a contradiction that the assertions of the
lemma are not true.

Set $L := V(G) - V(F)$ and $m := |L|$. Let $\mathcal{C} :=
\{V_1(K):K \in F\}, \mathcal{D} := \{V_2(K): K \in F\}$
and $C := \bigcup \mathcal{C}, D:=\bigcup \mathcal{D}$.
We call members of $\mathcal C$ \emph{lilliputs} while members of $\mathcal D$ are
\emph{giants}. We say that giant $V_2(K)$ $(K\in F)$ is
\emph{coupled} with lilliput $V_1(K)$.

As $F$ is a maximum \rcolor{size} $K_{s,t}$-tiling in $G$,
by~\eqref{eq:KTS} we have that 
\begin{equation}\label{eq:H0}
e(G[L]) = o(n^2)\;.
\end{equation}

Let $r$ be the number of copies of $K_{s,t}$ in $F$. Then
$r \leq (1-\eps)\alpha n / (s+t)$. Moreover, we have
\begin{equation}\label{eq:DefM}
m=n-(s+t)r\;.
\end{equation}

Let us define an auxiliary graph $H=(V',E')$ as follows.
The vertex-set of $H$ is $V' := \mathcal{C} \cup
\mathcal{D} \cup L$. For any $x \in L$ and $K \in F$ the edge $xV_1(K)$
belongs to $E'$ iff $\neighbor_G(x) \cap V_1(K) \neq
\emptyset$. Similarly, the edge $xV_2(K)$ belongs to $E'$
iff $\neighbor_G(x) \cap V_2(K) \neq \emptyset$. Finally,
for any distinct $K,K' \in F$ the edge $V_2(K)V_2(K')$
belongs to $E'$ iff $E_G(V_2(K),V_2(K')) \neq \emptyset$. The vertices $L$
and the vertices $\mathcal C$ induce two independent sets
in $H$.

As~(i) does not hold, $H[L,\mathcal{D}]$ does not contain
a matching with at least $\eps'n$ edges. It follows from
Fact~\ref{bip_matching} that 
\begin{equation}\label{eq:H1}
e_G(L, D) \leq \eps'nt\max\{m,r\} \rcolor{\leq t\eps'n^2}\;.
\end{equation}

Let $M$ be a maximum matching in $H[L, \mathcal{C}]$ with
$l$ edges. Obviously, $l\le r$. By
Fact~\ref{bip_matching}, we have that 
\begin{equation}\label{eq:H2}
e_G(L, C) \leq ls\max\{m,r\}\;.
\end{equation} Let
$\mathcal C'\subseteq \mathcal C$ be the lilliputs matched by $M$. We write $\mathcal D'\subseteq\mathcal D$ for the giants coupled with $\mathcal C'$. Set $D'=\bigcup
\mathcal D'$.

Suppose for a moment that $H[\mathcal{D}'] \cup
H[\mathcal{D}',\mathcal{D}-\mathcal{D}']$ contains a matching $T$
with at least $\eps'n$ edges. Let $\mathcal{D}''$ be the giants \rcolor{in $\mathcal{D}'$} 
matched by $T$ and $M'$ the set of edges in $M$ matching the 
lilliputs coupled with $\mathcal{D}''$. Then $M'$ and $T$ give rise 
to an $F$-augmentation $\rcolor{(}E_0, E_1\rcolor{)}$ in $G$ with $|E_0| = |M'| \geq |T| \geq \eps'n$, contradicting our assumption that (ii) does not hold.

So $H[\mathcal{D}'] \cup
H[\mathcal{D}',\mathcal{D}-\mathcal{D}']$ does not
contain a matching with at least $\eps'n$ edges. Applying Theorem~\ref{thm:ErdGall} and passing to the graph $G$,
we get
\begin{equation*}
e(G[D'] \cup G[D', D - D'])\le
t^2\ex(r,\eps'n\rcolor{ \times K_2})+r{t\choose 2}\le
2t^2\eps'nr+r{t\choose 2}\;.
\end{equation*}
Therefore,
\begin{align}\nonumber
e(G[C\cup D])&=e(G[D'] \cup G[D', D - D'])+e(G[D-D'])+e(G[C])+e_G(C,D)\\ 
\label{eq:H3}
&\le 2t^2\eps'nr+r{t\choose 2}+\rcolor{{(r-l)t\choose 2}}+{rs\choose 2}+r^2st.
\end{align}
Summing up the
bounds~\eqref{eq:H0},~\eqref{eq:H1},~\eqref{eq:H2},
and~\eqref{eq:H3} we get:
\begin{alignat*}{1}
	e(G) & = e(G[L]) + e_G(L,D) + e_G(L,C) + e(G[C \cup D]) \\
	& \leq o(n^2) + t\eps'n^2 + ls\max\{m,r\} + 2\eps'nrt^2
	+ r{t \choose 2}  + \rcolor{{(r-l)t\choose 2}} + {rs \choose 2} + r^2st. \\
\intertext{\rcolor{Using the convexity of $f(l) := ls\max\{m, r\}+ {(r-l)t\choose 2}$ on $[0, r]$, and the fact that $rt \leq n$, 
we get:}}
	e(G) &\leq o(n^2) + 3t\eps'n^2 + r{t \choose 2}
	+ r^2st + {rs \choose 2} + \max\left\{{rt \choose 2}, rs\max\{m, r\}\right\}. \\
\intertext{\rcolor{However, $r^2s \leq {rt \choose 2} + o(n^2)$, and hence from \eqref{eq:DefM} we get:}}
	e(G) &\leq o(n^2) + 3t\eps'n^2 + r{t \choose 2}
	+ r^2st + {rs \choose 2} + \max\left\{{rt \choose 2}, rs(n-(s+t)r)\right\} \\ 
	& < \max\left\{{(s+t)r \choose 2}, {rs \choose 2} + rs(n-rs)\right\} + (3t+1)\eps'n^2,
\end{alignat*}
\rcolor{
where in the last inequality we have majorized the term $r\binom{t}{2}+o(n^2)$ by $\epsilon' n^2$. But
\begin{equation*}
{(s+t)r \choose 2} + (3t+1)\eps'n^2 
\leq \binom{(1-\eps)\alpha n}{2} + \frac{\eps \alpha^2n^2}{4} 
< \left(1-\frac{\eps}{2}\right)\frac{\alpha^2n^2}{2},
\end{equation*}
and
\begin{alignat*}{1}
{rs \choose 2} + rs(n-rs) + (3t+1)\eps'n^2 
&< rsn - \frac{r^2s^2}{2}+\frac{\eps s\alpha n^2}{4(s+t)} \\
&\leq\frac{2s\alpha}{s+t}\left(1-\frac{\alpha s}{2(s+t)} +\frac{\eps(2-\eps)\alpha s}{2(s+t)}-\frac{3\eps}{4}\right)\frac{n^2}{2} \\
&< \frac{2s\alpha}{s+t}\left(1-\frac{\alpha s}{2(s+t)} -\frac{\eps}{4}\right)\frac{n^2}{2}.
\end{alignat*}
Consequently for large enough $n$,
}
\begin{equation*}
	e(G) < T_{s,t}(\alpha){n \choose 2}\;, 
\end{equation*}
a contradiction.
\end{proof}

Suppose $G=(V,E)$ is a graph and $r \in \mathbb N$.  The
\emph{$r$-expansion} of $G$ is the graph $G'=(V',E')$ defined as follows. The vertex set of $G'$ is $V \times [r]$. For $a,b \in [r]$, an edge $((u,a),(v,b))$ belongs to $E'$ iff $uv$ belongs to $E$. Note that there is a natural projection $\pi_{G'}: V' \rightarrow V$ that maps every vertex $(u,a)$ from $G'$ to the vertex $u$ in $G$. We are interested in the following property of $r$-expansions. Suppose that $K$ is a copy of any graph from $\mathcal
F^*$ in $G$. Then
$\pi_{G'}^{-1}(V(K))$ contains a complete bipartite graph
$B$ with color classes of sizes $s(K)r$ and $t(K)r$. By
Lemma~\ref{lem:til} we can tile $B$ almost perfectly with
copies of $K_{s,t}$. If $F$ is an $\mathcal F^*$-tiling
in $G$, we can apply the above operation on
each member $K\in F$ and obtain a new tiling $F'$ ---
which we call \emph{retiling} --- in the graph
$G'$.

We are now ready to prove Theorem~\ref{main_result}.


\begin{proof}[Proof of Theorem~\ref{main_result}]
Note that it suffices to prove the theorem for $H \simeq K_{s,t}$.

We first deal with the particular case $t =
s$. Set $\alpha':=(1-\epsilon/4)\alpha$. Let
$\eps_1:=\frac1{\rcolor{5}}(T_{s,t}(\alpha)-T_{s,t}(\alpha'))$,
and $\eps_2$ be given by Lemma~\ref{lem:tiling2} for input
parameters $H$, $d:=\eps_1$ and
$\gamma:=\alpha\epsilon/8$. Suppose that $k_0$ is
sufficiently large. Let $M$ be the bound from
Lemma~\ref{lem:reg} for precision $\eps_R:=\min\{\eps_1,\eps_2\}$ and minimal number of
clusters $k_0$. Let $C$ be given by Lemma~\ref{lem:til} for
the input parameters $s,t$. Fix $n_0\gg MC$. Suppose that
$G$ is an $n$-vertex graph, $n\ge n_0$, with at least
$T_{s,t}(\alpha){n\choose 2}$ edges. We apply Lemma~\ref{lem:reg}
on $G$ to obtain an $(\eps_R,d)$-reduced graph $R$ with
$k$ clusters, $k_0\le k\le M$. We have that
$$e(R)\ge (T_{s,t}(\alpha)-d-3\eps_1){k\choose
2}=(T_{s,t}(\alpha')+\frac1{\rcolor{5}}(T_{s,t}(\alpha)-(T_{s,t}(\alpha'))){k\choose
2}\gByRef{eq:TssVersusErdosGallai}\ex\left(k,\frac{\alpha'
k}2\times K_2\right)\;.$$
Therefore, $R$ contains at least $\frac{\alpha'
k}2$ independent edges. These edges correspond to regular
pairs in $G$ which can be tiled
almost perfectly with copies of $K_{s,t}$, by means of  Lemma~\ref{lem:til} and Lemma~\ref{lem:tiling2}. Elementary
calculations give that in this way we get a tiling of size
at least $(1- \eps)\alpha n$.
%

Consequently we may suppose that $t > s$.
We first define a handful of parameters.
Set 
\begin{equation*}
\label{eq:param1}
	\alpha':= \frac{6-4\eps}{6-3\eps}\alpha,
	\quad \gamma := (1-\eps/2)\alpha', 
	\quad d := \frac25(T_{s,t}(\alpha)-T_{s,t}(\alpha'))\;.
\end{equation*}
Note that $\gamma = (1 - 2\eps/3)\alpha$.

Let
$\eps_R$ be given by Lemma~\ref{lem:key} for input
graph $K_{s,t}$, density $d/2$ and approximation parameter $\gamma$. We may suppose that $\eps_R$ is sufficiently small such that $\gamma(1-\eps_R) > (1-\eps)\alpha$ and $\eps_R < d/2$. Let $C$ be given by Lemma~\ref{lem:til} for input $s,t$. Further, let $\eps'$ and $h$ be given by
Lemma~\ref{main_lemma} for input parameters $\alpha'$ and
$\eps/4$. \rcolor{We may assume that $\eps' < \eps$.}
Set
\begin{equation*}
\label{eq:param2}
	p := t^2\left\lceil\frac{4C}{\eps'}\right\rceil, 
	\quad q := \left\lceil\frac{2t}{\eps'}\right\rceil
\end{equation*}

Let $M$ be the
upper bound on the number of clusters given by
Lemma~\ref{lem:reg} for input parameters $h$ (for the
minimal number of clusters) and $\eps_Rp^{-q}/2$ (for
the precision). Let $n_0 > Mp^q$ be sufficiently large.

Suppose now that $G$ is a graph with $n>n_0$ vertices and
at least $T_{s,t}(\alpha){n\choose 2}$ edges. We first apply
Lemma~\ref{lem:reg} to $G$ with parameters $\eps_Rp^{-q}/2$ and $h$. In this way we obtain an $(\eps_Rp^{-q}/2,d)$-reduced graph $R$ with at least $h$ vertices. 

Let us now define a sequence of graphs
$R^{(i)}$ by setting $R^{(0)} = R$ and letting $R^{(i)}$ be the $p$-expansion of $R^{(i-1)}, i = 1,2,\ldots,q$. 
Note that
$e(R^{(i)})\ge T_{s,t}(\alpha'){\order{R^{(i)}}\choose 2}$ for
every $i\in\{0,1,\ldots,q\}$.

Let $F^{(i)}$ be a maximum \rcolor{size} $K_{s,t}$-tiling in $R^{(i)}$
for $i=0,1,\ldots,q$. We claim that
\begin{equation}
\label{eq:Improving}
|F^{(i)}| \ge 
	\min\left\{
	\frac{i\eps' \order{R^{(i)}}}{2t},
	\left(1-\frac\eps2\right)\alpha' \order{R^{(i)}}
	\right\}\;.
\end{equation}
To this end it suffices to show that for any $i \geq 1$,
\begin{enumerate}
\item[(C1)] if 
$|F^{(i-1)}| > (1-\eps/4)\alpha'\order{R^{(i-1)}}$,
then
$\frac{|F^{(i)}|}{\order{R^{(i)}}} \ge
\frac{|F^{(i-1)}|}{\order{R^{(i-1)}}}-\frac{\eps\alpha'}4$, and
\item[(C2)]
if $|F^{(i-1)}|\le
(1-\eps/4)\alpha' \order{R^{(i-1)}}$, then
$\frac{|F^{(i)}|}{\order{R^{(i)}}}\ge
\frac{|F^{(i-1)}|}{\order{R^{(i-1)}}}+\frac{\eps'}{2t}$.
\end{enumerate}

In the case (C1), \rcolor{according to Lemma \ref{lem:til}}, the retiling of $F^{(i-1)}$ in $R^{(i)}$ has size 
at least $\rcolor{|F^{(i-1)}|(p-C)>}(1-\eps/2)\alpha' \order{R^{(i)}}$, thus proving the
statement.

Consequently we may suppose that we are in case (C2). Apply 
Lemma~\ref{main_lemma} to the graph $R^{(i-1)}$ and the tiling 
$F^{(i-1)}$, with parameters $\alpha'$ and $\eps/4$.

Suppose first that assertion~(i) of the lemma holds. Then 
$R^{(i-1)}$ contains an
$\mathcal{F}_1$-tiling $F$ with $\frac{|F|}{\order{R^{(i-1)}}}
\geq \frac{|F^{(i-1)}|}{\order{R^{(i-1)}}} + \eps'$.  By
retiling $F$, we get a $K_{s,t}$-tiling in $R^{(i)}$ with size at least $|F|(p-C) > 
i\eps'\order{R^{(i)}}/(2t)$,
thus proving the statement.

Suppose now that assertion~(ii)~of Lemma~\ref{main_lemma} is true.
Then $R^{(i-1)}$ contains an $F^{(i-1)}$-augmentation
$\rcolor{(}E_0, E_1\rcolor{)}$ with $|E_0| \geq \eps'\order{R^{(i-1)}}$. Let $r=p/t$. We shall denote by $T$ the $t$-expansion of $R^{(i-1)}$ and by $T'$ the $r$-expansion 
of $T$. Note that $T'$ is isomorphic to $R^{(i)}$.

Let us build an $\mathcal{F}_2$-tiling in $T$ in the following way. 

For every edge $e = (u,v) \in E_0$ with $u \in V(F^{(i-1)})$ 
we choose an edge $e'=(u',v')$ in $T$ with $\pi_{T}(u') = u$ and 
$\pi_{T}(v') = v$ . We shall denote by $w_e$ the vertex $u'$ 
corresponding to $u$. 

For every edge $e = (u,v) \in E_1$ we choose a set $S_e$ of $t$ independent edges in $\pi_{T}^{-1}(e)$.

For every $K \in F^{(i-1)}$ we shall also choose a 
subgraph $K'$ of $T$. We  distinguish the following
cases. If $K$ has no vertex  matched by $E_0$ or $E_1$,
then we let $K' := T[\pi_{T}^{-1}(K)]$. If $K$  has a
vertex $u$ matched by $E_1$ but no vertex matched  by
$E_0$, we let $K': = T[\pi_{T}^{-1}(K-u)]$. Then  $K'
\simeq K_{st,(t-1)t}$. Finally, if $K$ has a vertex  $u$
matched by an edge $e \in E_0$ and a vertex $v$ matched 
by an edge in $E_1$, we let $K' := T[\pi_{T}^{-1}(K-v)] -
w_e$.  Note that in this last case $K' \simeq K_{st-1,(t-1)t}$.

It is easy to see that 
\begin{equation*}
\label{eq:tiling}
F := \{e': e \in E_0\} \cup \{K' : K \in F^{(i-1)}\} \cup
\left(\bigcup_{e \in E_1}S_e\right)
\end{equation*}
is an $\mathcal{F}_2$-tiling in $T$. Moreover, we have that 
$\frac{|F|}{\order{T}} \geq \frac{|F^{(i-1)}|}{\order{R^{(i-1)}}} +
\frac{\eps'}{t}$. So the retiling of $F$ in $T'$ has size
at least $|F|(r-C) \geq i\eps'\order{R^{(i)}}/(2t)$. This proves (C2) and also~\eqref{eq:Improving}.

Using Lemma~\ref{lem:slicing}, we may subdivide every cluster 
corresponding to a vertex of $R$ into $p^q$ equal-sized parts, by 
discarding some vertices if necessary. This gives us an 
$(\eps_R,d/2)$-reduced graph $R'$. By construction $R' \simeq 
R^{(q)}$. By~\eqref{eq:Improving}, there is a $K_{s,t}$-tiling 
$F$ in $R'$ with size at least $(1-\eps/2)\alpha'\order{R'}$.
Let $G'$ be the subgraph of $G$ induced by the clusters corresponding to 
the vertices of $R'$. By applying Lemma~\ref{lem:key} to $R'$, we see 
that $G'$ has a $K_{s,t}$-tiling of size at least $\gamma
\order{G'} \geq \gamma(1-\eps_R)\order{G} > (1-\eps)\alpha \order{G}$,
and so does $G$.

This finishes the proof of Theorem~\ref{main_result}.
\end{proof}
\section*{Acknowledgements}
We would like to thank the anonymous referee for his/her very useful comments.

\bibliographystyle{amsplain} \bibliography{bibl}
\end{document}